

%

%

%


\documentstyle{amsppt}
\magnification=\magstep1
\NoBlackBoxes


\topmatter

\title A correlation inequality  
for the expectations of norms of stable vectors.
\endtitle

\author Alexander Koldobsky \endauthor 
\address Division of Mathematics and Statistics, 
University
of Texas at San Antonio, San Antonio, TX 78249, U.S.A. \endaddress
\email koldobsk\@ringer.cs.utsa.edu \endemail

\rightheadtext{A correlation inequality}

\abstract 
For $0<q\le 2,\ 1\le k < n,$ let 
$X=(X_1,...,X_n)$ and $Y=(Y_1,...,Y_n)$ be symmetric 
$q$-stable random vectors so that the joint distributions of 
$X_1,...,X_k$ and  $X_{k+1},...,X_n$ are equal to the joint 
distributions of $Y_1,...,Y_k$ and $Y_{k+1},...,Y_n,$ respectively, 
but $Y_i$ and $Y_j$ are independent for every 
$1\le i \le k,\ k+1\le j \le n.$ We prove that
$\Bbb E (f(X)) \ge \Bbb E (f(Y))$
where $f$ is any continuous,
positive, homogeneous of the order $p\in (-n,0)$ function on 
$\Bbb R^n\setminus \{0\}$
such that $f$ is a positive definite distribution in $\Bbb R^n,$ 
and $f(u,v)=f(u,-v)$ 
for every $u\in \Bbb R^k,\ v\in \Bbb R^{n-k}.$
As a particular case, we show that
$$\Bbb E\ (\max_{i=1,...,n} |X_i|)^p \ge 
\Bbb E\ (\max_{i=1,...,n} |Y_i|)^p$$
for every $p\in (-n,-n+1).$ The latter inequality is
related to Slepian's Lemma and to the Gaussian correlation
problem.
\endabstract


\rightheadtext{An inequality for stable expectations}

\thanks Part of this work was done during
the NSF Workshop in Linear Analysis and Probability held 
at Texas A\&M University in Summer 1995 \endthanks

\endtopmatter \document \baselineskip=14pt

\head 1. Introduction \endhead
Let $X_1,...,X_n$ and $Y_1,...,Y_n$ be symmetric jointly Gaussian 
random variables. The well-known Slepian's Lemma  \cite{20}, 
\cite{9} states that
if $\Bbb E X_i^2 = \Bbb E Y_i^2$ and $\Bbb E (X_iX_j) \ge \Bbb E (Y_iY_j)$
for every $i,j=1,...,n,$ then 
$$\Bbb E\ (\max_{i=1,...,n} X_i) \le \Bbb E\ (\max_{i=1,...,n} Y_i), $$
and, even more, for every $t\in \Bbb R$ one has
$$P(\max_{i=1,...,n} X_i > t) \le P(\max_{i=1,...,n} Y_i > t),$$
These inequalities mean that the maximum of Gaussians tends to be larger 
when they are less correlated.
\medbreak
One may ask a question of whether the absolute values of Gaussians
behave in a similar way, namely, given a number $k\in N,\ 1\le k < n$ 
and fixed 
joint distributions of $X_1,...,X_k$ and of $X_{k+1},...,X_n$ is it true
that the quantities $\Bbb E\ (\max_{i=1,...,n} |X_i|)$ and 
$P(\max_{i=1,...,n} |X_i| > t)$ are maximal when the random variables
$X_i$ and $X_j$ are independent for every $1\le i \le k,\ k+1\le j \le n?$
However, this question for $P(\max_{i=1,...,n} |X_i| > t)$ appears to 
be equivalent 
to the famous correlation  problem for Gaussian measures of symmetric 
convex sets: Is it true that $\mu(A\cap B)\ge \mu(A) \mu(B)$ for any 
symmetric convex sets $A, B$ in
$\Bbb R^d$ where $\mu$ is the standard symmetric Gaussian measure in 
$\Bbb R^d ?$ Pitt \cite{16} has proved that the answer to the correlation
problem is positive in the case $d=2,$ and, therefore, confirmed the proper
behaviour of the Gaussian random variables generated by a two-dimensional 
Gaussian vector. For $d>2$ the problem remains open 
(see \cite{19} for
the history of the problem and partial results).

\bigbreak

This note provides some evidence supporting the conjecture on the 
behaviour of the absolute values of Gaussians. 
Let $0<q\le 2,\ 1\le k < n.$ Throughout the paper,
$X=(X_1,...,X_n)$ and $Y=(Y_1,...,Y_n)$ are symmetric 
$q$-stable random vectors such that the joint distributions of 
$X_1,...,X_k$ and  $X_{k+1},...,X_n$ are equal to the joint 
distributions of $Y_1,...,Y_k$ and $Y_{k+1},...,Y_n,$ respectively, 
but $Y_i$ and $Y_j$ are independent for every 
$1\le i \le k,\ k+1\le j \le n.$ We prove that in this situation
$$\Bbb E\ (\max_{i=1,...,n} |X_i|)^p \ge 
\Bbb E\ (\max_{i=1,...,n} |Y_i|)^p, \tag{1}$$
for every $p\in (-n,-n+1).$ 
\medbreak
We show this result as a particular case of the following more general 
inequality: 
\proclaim{Theorem 1} Let $q,k, X,Y$ be as above,
and let $-n<p<0$ and $f$ be a continuous,
positive, homogeneous of the order $p$ function on 
$\Bbb R^n\setminus \{0\}$
such that $f$ is a positive definite distribution in $\Bbb R^n,$ 
and $f(u,v)=f(u,-v)$ 
for every $u\in \Bbb R^k,\ v\in \Bbb R^{n-k}.$ Then
$\Bbb E (f(X)) \ge \Bbb E (f(Y)).$ \endproclaim

\medbreak

The inequality (1) will follow from Theorem 1 and a simple 
fact that, for $-n<p<-n+1,$  every positive, continuous, 
homogeneous of the order $p$ function on $\Bbb R^n\setminus \{0\}$ 
is also a positive definite distribution.

\bigbreak

We refer the reader to \cite{7, 17, 18} for other
results related to the Slepian Lemma.
\bigbreak

\head 2. Expectations of positive powers of norms \endhead

We start with an inequality for the expectations of positive powers
of certain norms. The techniques used in this case are 
standard, but the positive case makes more clear what happens 
later in the case of negative powers.

We need a few simple inequalities for the $L_q$-norms which follow
from Clarkson's inequality (see \cite{2}). For the reader's 
convenience
we include the proof. We denote by $\|\cdot\|_q$ the norm of the 
space $L_q([0,1]).$
\bigbreak
\proclaim{Lemma 1} Let $x,y\in L_q([0,1]),\ 0<q\le 2.$ Then
$$\exp(-\|x+y\|_q^q) + \exp(-\|x-y\|_q^q) \ge
2\exp(-\|x\|_q^q-\|y\|_q^q). \tag{2}$$
Also for every $0<p\le q$
$$\|x + y\|_q^p +
\|x - y\|_q^p \le
2(\|x\|_q^q + \|y\|_q^q)^{p/q}.\tag{3}$$
Finally, for $q=2$ and $p>2$ the inequality (3) goes in the 
opposite direction.
\endproclaim

\demo{Proof} First, note that for any $0<q\le 2$
$$\|x+y\|_q^q + \|x-y\|_q^q \le 2(\|x\|_q^q + \|y\|_q^q), \tag{4}$$
and this is a simple consequence of the same inequality for real
numbers. Now to get (2) apply the relation between the
arithmetic and geometric means and then use (4). The inequality (3)
also follows from (4):
$$ \Big( {{\|x + y\|_q^p + \|x - y\|_q^p}\over 2} \Big)^{1/p} \le
\Big( {{\|x + y\|_q^q + \|x - y\|_q^q}\over 2} \Big)^{1/q} \le
(\|x\|_q^q + \|y\|_q^q)^{1/q}.$$
Finally, if $q=2$ the latter calculation works for $p>2$ where
the first inequality goes in the opposite direction, and the 
second inequality turns into an equality.
\qed \enddemo

\bigbreak

For $0<q\le 2,\ 1\le k < n,$ let $X=(X_1,...,X_n),\ Y=(Y_1,...,Y_n)$ 
be the symmetric $q$-stable random vectors defined in Introduction.

The characteristic function of the vector $X$ has the form
$$\phi(\xi)= \exp(-\|\sum_{i=1}^n \xi_is_i\|_q^q),\quad 
\xi\in \Bbb R^n, \tag{5}$$
where $s_1,\dots,s_n\in L_q([0,1]).$

Then the characteristic function of $Y$ is equal to
$$\phi_0(\xi)= \exp(-\|\sum_{i=1}^k \xi_is_i\|_q^q - 
\|\sum_{i=k+1}^n \xi_is_i\|_q^q).$$

\medbreak

Let $(\Bbb R^n, \|\cdot\|)$ be an $n$-dimensional subspace of 
$L_p([0,1]),\ p>0.$
A well-known easy fact due to P.Levy \cite{13} is that an
$n$-dimensional space is isometric to a subspace of $L_p([0,1])$ 
if and only if its norm admits the following Levy representation:
$$\|x\|^p =\int_S |(x,\xi)|^p\ d\gamma(\xi) \tag{6}$$
for every $x\in \Bbb R^n,$ where $S$ is the unit sphere in $\Bbb R^n,$
$(x,\xi)$ stands for the scalar product, 
and $\gamma$ is a finite Borel (non-negative) measure on $S.$

\proclaim{Proposition 1} Let $q, k, X, Y$ be as in the Introduction,
$0<p\le q$ and $(\Bbb R^n,\|\cdot\|)$ is a 
subspace of $L_p$ with the norm satisfying $\|(u,v)\|= \|(u,-v)\|$
for every $u\in \Bbb R^k,\ v\in \Bbb R^{n-k}.$ 
Then $\Bbb E\ (\|X\|^p) \le \Bbb E\ (\|Y\|^p).$ 
Also if $q=2$ and $p>2$ the inequality goes in the opposite direction. 
\endproclaim
 
\demo{Proof} A basic property of the stable vector with the characteristic 
function (5) is that, for any
vector $\xi\in \Bbb R^n,$ the random variable $(X,\xi)$ has the same
distribution as $\|\sum_{i=1}^n \xi_is_i\|_q Z,$ where $Z$ is the
standard one-dimensional $q$-stable random variable.
Therefore, if $p<q$  then 
$$\Bbb E\ |(X,\xi)|^p = c_{p,q}\|\sum_{i=1}^n \xi_is_i\|_q^p,\tag{7}$$
where $c_{p,q}$ is the $p$-th moment of $Z$
(which exists only for $p<q$ if $q<2,$ and
it exists for every $p>0$ if $q=2;$ see \cite {22} for
a formula for $c_{p,q}).$ Similarly, we get
$$\Bbb E\ |(X_{-},\xi)|^p = c_{p,q}\|\sum_{i=1}^k \xi_is_i- 
\sum_{i=k+1}^n \xi_is_i\|_q^p, $$ 
where $X_{-}=(X_1,...,X_k,-X_{k+1},...,-X_n).$ Also,
$$\Bbb E\ |(Y,\xi)|^p = c_{p,q}(\|\sum_{i=1}^k \xi_is_i\|_q^q + 
\|\sum_{i=k+1}^n \xi_is_i\|_q^q)^{p/q}.$$

Since $(\Bbb R^n,\|\cdot\|)$ is a subspace of $L_p([0,1]),$
we can use the Levy representation (6) and after that 
the formula (7) to get

$$\Bbb E\ (\|X\|^p) = 
\int_S \Bbb E\ (|(X,\xi)|^p)\ d\gamma(\xi)=
c_{p,q}\int_S \|\sum_{i=1}^n \xi_is_i\|_q^p\ d\gamma(\xi).\tag{8}$$
Similarly,
$$E (\|Y\|^p)= c_{p,q}\int_S (\|\sum_{i=1}^k \xi_is_i\|_q^q + 
\|\sum_{i=k+1}^n \xi_is_i\|_q^q)^{p/q}\ d\gamma(\xi),\tag{9}$$
$$\Bbb E\ (\|X_{-}\|^p)= c_{p,q}\int_S \|\sum_{i=1}^k \xi_is_i- 
\sum_{i=k+1}^n \xi_is_i\|_q^p \ d\gamma(\xi).\tag{10}$$
Since $0<p\le q,$ the equalities (8), (9), (10) in conjunction 
with (3) imply 
$\Bbb E\ (\|X\|^p)+\Bbb E\ (\|X_{-}\|^p)\le 2\Bbb E\ (\|Y\|^p),$
and now the result follows from the property of
the norm that $\|X\|=\|X_{-}\|.$ In the case $q=2,\ p>2$ we use
the corresponding part of Lemma 1.
\qed \enddemo

\bigbreak

\subheading{Remarks} (i) For $p>q,\ q< 2 $ the expectation
of $\|X\|^p$ does not exist so the statement of Proposition 1
does not make sense in that case.

\medbreak

(ii) In view of Proposition 1, it is natural to ask
how can one check whether a given space is isometric to a subspace 
of $L_p.$ This question is the matter of an old problem raised by 
P.Levy \cite{13}. 
In the same paper P.Levy showed that an $n$-dimensional 
space is isometric to a subspace of $L_p$ if and only if its norm 
admits the representation (6). Since then a few criteria involving 
the Fourier transform have appeared. Bretagnolle, Dacunha-Castelle
and Krivine \cite{1} proved that,
for $0<p\le 2,$ a Banach space is isometric to a subspace of $L_p$
if and only if the function $\exp(-\|x\|^p)$ is positive definite,
and, in particular, showed that the spaces $L_q$ embed isometrically
into $L_p$ if $0<p<q\le 2.$ Another Fourier transform criterion was given
in \cite{10}, \cite{11}: for any 
$p\in (0,\infty)\setminus \{even\ integers\},$
an $n$-dimensional space is isometric to a subspace of $L_p$ if and
only if the restriction of the Fourier transform of $\|x\|^p \Gamma(-p/2)$ 
to the unit sphere $S$ in $\Bbb R^n$ is a finite Borel measure on $S$
(the Fourier transform is considered in the sense of distributions).
Recently, two criteria were shown that were in terms of the derivatives
of the norm: Zastavny \cite{21} proved that a three dimensional space 
is not isometric to a subspace of $L_p$ with $0<p\le 2$ if there exists 
a basis $e_1,e_2,e_3$ so that the function
$$(y,z)\mapsto \|xe_1 + ye_2 + ze_3\|^{'}_x (1,y,z)/\|e_1+ye_2+ze_3\|,
\ y,z \in \Bbb R$$
belongs to the space $L_1(\Bbb R^2).$ 
By inverting the representation (6), it was shown in \cite{12} that
an $n$-dimensional space is isometric to a subspace of $L_p$ with 
$n+[p]$ being an even integer if
$(-1)^{(n+[p])/2}\Delta^{(n+[p])/2} \|x\|^p$ is a positive continuous 
function on the sphere $S$ ( $\Delta$ is the Laplace operator, $p$ is
not an even integer; in the case 
where $n+[p]$ is odd the formula must be slightly modified.)

\medbreak

(iii) Misiewicz \cite{15} proved that the spaces $\ell_{\infty}^n,\ n>2$
do not embed in any of the spaces $L_p, \ p>0,$ therefore Proposition 1  
does not tell anything about the behaviour of $\max(|X_i|)$ (except
for the case $n=2$ where one can use the well-known fact due to
Herz \cite{8}, Ferguson \cite{3}, Lindenstrauss \cite{14} that
any two-dimensional Banach space embeds 
isometrically in each one of the spaces $L_p$ with $0<p\le 1.$)

\head 3. Expectations of negative powers of norms \endhead

As one can see from Remarks (ii) and (iii) the condition
of Proposition 1 that the norm embeds isometrically in $L_p$
is quite restricting and is not easy to check. In this section
we replace this condition by an equivalent one, and that allows
us to extend the result of Proposition 1 to the case of negative
powers $p$ and, more important, to a much larger class of norms.

\medbreak

To formulate the equivalent condition, we need some notation.
As usual, we denote by $\Cal S(\Bbb R^n)$ the space of rapidly 
decreasing infinitely
differentiable functions in $\Bbb R^n,$ and by $\Cal S^{'}(\Bbb R^n)$
the space of distributions over $\Cal S(\Bbb R^n).$ Recall that 
the Fourier transform of a distribution $f\in \Cal S^{'}(\Bbb R^n)$ 
is defined by $(\hat{f},\phi) = (f,\hat{\phi})$ for every
$\phi\in \Cal S(\Bbb R^n).$ We say that a 
distribution $f \in \Cal S^{'}(\Bbb R^n)$ 
is positive definite in a domain $D\subset \Bbb R^n$ if 
the Fourier transform of $f$ is a positive distribution
in $D,$ namely, $(\hat{f},\phi)\ge 0$ for every non-negative
function $\phi\in \Cal S(\Bbb R^n)$ supported in $D.$ 
We need the following

\proclaim{Lemma 2} 
Let $p\in (-1,\infty)$, $p$ is not an even integer. Let $\phi$ be a 
function from the space $\Cal S(\Bbb R^n).$  Then for every
$\xi\in\Bbb R^n$, $\xi\not = 0$,
$$
\int_{\Bbb R^n} |(x,\xi)|^p \hat{\phi}(x)\ dx=
(2\pi)^{n-1}c_p
(|z|^{-1-p}, \phi(z\xi))$$
where $c_p=(2^{p+1}\pi^{1/2}\Gamma ((p+1)/2))/{\Gamma(-p/2)},$
and $(|z|^{-1-p},\phi(z\xi))$ is the value 
of the one-dimensional distribution $|z|^{-1-p}$ at the test function
$z\to \phi(z\xi),\ z\in \Bbb R.$ 
\endproclaim

\demo{Proof} By the Fubini theorem
$$
\int_{\Bbb R^n}|( x,\xi)|^p \hat\phi(x)\, dx=
\int_{\Bbb R}|t|^p
\Big(\int_{(x,\xi)=t} \hat \phi (x)\, dx\Big)\, dz=
\big( |t|^p,
\int_{(x, \xi) =t} \hat \phi (x)\, dx\big).
\tag{11}
$$
It is well--known that the Fourier transform of
the function $t\to |t|^p,\ t\in \Bbb R$ is equal to
$(|t|^p)^\wedge (z)=c_p|z|^{-1-p}$,
$z\not= 0$, for every $p\in (-1,\infty)$ which is not an even integer
(see \cite{4}). Also the function $z\to (2\pi)^n \phi(-z\xi)$ 
is the Fourier transform of the function 
$$t\to \int_{(x, \xi)=t} \hat\phi (x)\, dx$$ 
(this is the connection between the Fourier transform
and the Radon transform, see \cite{6}).
Passing to the Fourier transforms in the equality (11) we get
$$ \big(|t|^p, \int_{(x, \xi) =t} \hat \phi (x)\, dx\big)=
(1/2\pi) \big( c_p|z|^{-1-p}, (2\pi)^n \phi(z\xi)\big). \qed
$$ \enddemo
\bigbreak
Now we are able to show the equivalent condition mentioned above:

\proclaim{Lemma 3} Let $p$ be a positive number which is 
not an even integer. A space $(\Bbb R^n,\|\cdot\|)$ is   
isometric to a subspace of $L_p([0,1])$ if and only if there
exists a finite Borel (non-negative) measure $\gamma$ on the unit sphere
$S$ in $\Bbb R^n$ so that, for every $\phi\in \Cal S(\Bbb R^n),$
$$\big((\|x\|^p)^{\wedge},\phi\big) =
c_p \int_S (|z|^{-1-p},\phi(z\xi))\ d\gamma(\xi).\tag{12}$$
\endproclaim

\demo{Proof} A simple fact going back to P.Levy \cite{13} is that
a space $(\Bbb R^n,\|\cdot\|)$ is isometric to a subspace of
$L_p([0,1])$ if and only if the norm admits the Levy representation
(6) with a measure $\gamma$ on the sphere $S.$ By (6) and Lemma 2, 
the space embeds into $L_p([0,1])$ if and only if, for every 
$\phi\in \Cal S(\Bbb R^n),$  
$$\big((\|x\|^p)^{\wedge},\phi\big) =
(\|x\|^p,\hat{\phi}) = \int_{\Bbb R^n} \|x\|^p \ \hat{\phi}(x)\ dx =$$
$$\int_S d\gamma(\xi)\ \Big(\int_{\Bbb R^n} 
|(x,\xi)|^p\ \hat{\phi}(x)\ dx \Big) = 
c_p \int_S (|z|^{-1-p},\phi(z\xi))\ d\gamma(\xi). \qed$$ 
\enddemo

\bigbreak

If the function $\phi$ in (12) is supported in 
$\Bbb R^n\setminus \{0\},$
we have $(|z|^{-1-p},\phi(z\xi))=\int_{\Bbb R} |z|^{-1-p}\phi(z\xi)$
which is non-negative if the function $\phi$ is non-negative. 
Therefore,

\proclaim{Corollary 1} If a space $(\Bbb R^n,\|\cdot\|)$ embeds 
isometrically in $L_p([0,1])$ with $p>0,\ p\neq 2k,\ k\in \Bbb N,$
then the distribution $\|x\|^p \Gamma(-p/2)$ is positive definite
in $\Bbb R^n\setminus \{0\}.$ 
\endproclaim

\bigbreak

Now we are able to prove an analog of Proposition 1 for 
negative powers $p$ replacing the embedding in $L_p$ by
positive definiteness of $\|x\|^p$ (note that for negative
$p$ the numbers $\Gamma(-p/2)$ are always positive). Also in the
case of negative $p$ we will be able to replace the norm to the
power $p$ by any positive, continuous, homogeneous of the order $p$ 
function. (The latter means that $f(tx)=|t|^p f(x)$ for every
$t\in \Bbb R,\ t\neq 0, x\in \Bbb R^n\setminus \{0\}.)$

\proclaim{Theorem 1} Let $q,k, X,Y$ be as in the Introduction,
and let $-n<p<0$ and $f$ be a continuous,
positive, homogeneous of the order $p$ function on 
$\Bbb R^n\setminus \{0\}$
such that $f$ is a positive definite distribution in $\Bbb R^n,$ 
and $f(u,v)=f(u,-v)$ 
for every $u\in \Bbb R^k,\ v\in \Bbb R^{n-k}.$ Then
$\Bbb E (f(X)) \ge \Bbb E (f(Y)).$ \endproclaim

\demo{Proof} We use the following generalization of Bochner's theorem
(see \cite{5}): if $f$ is a positive definite 
distribution in $\Bbb R^n$ (over 
$\Cal S(\Bbb R^n))$ then $f$ is the Fourier transform 
(in the sense of distributions) of a tempered measure $\mu$ in
$\Bbb R^n.$ (Recall that a measure is called tempered if 
$\int_{\Bbb R^n} (1+\|x\|_2)^{\alpha}\ d\mu(x) < \infty$ for
some $\alpha < 0.)$ Let $\mu$ be the tempered measure whose 
Fourier transform is equal to $f.$ 

Let $P_X$ be the $q$-stable measure in $\Bbb R^n$ 
according to which the random vector $X$ is distributed.
Applying the Parseval equality and the expression (5) for the 
characteristic function of $X$ we get
$$\Bbb E (f(X))= \int_{\Bbb R^n} f(x)\ dP_X(x) =
\int_{\Bbb R^n} \widehat{P_X}(\xi)\ d\mu(\xi) =
\int_{\Bbb R^n} \ exp(-\|\sum_{i=1}^n \xi_is_i\|_q^q)\ d\mu(\xi).$$
Note that the function $f$ is locally integrable in $\Bbb R^n$ 
because $-n<p<0.$ Similarly,
$$\Bbb E (f(X_{-}))= 
\int_{\Bbb R^n} \ exp(-\|\sum_{i=1}^k \xi_is_i-
\sum_{i=k+1}^n \xi_is_i\|_q^q)\ d\mu(\xi),$$
where $X_{-}=(X_1,...,X_k,-X_{k+1},...,-X_n),$ and 
$$\Bbb E (f(Y))= 
\int_{\Bbb R^n} \ exp(-\|\sum_{i=1}^k \xi_is_i\|_q^q-
\|\sum_{i=k+1}^n \xi_is_i\|_q^q)\ d\mu(\xi).$$
Now by the inequality (2) from Lemma 1 and taking in account that
$\mu$ is a positive measure, we get
$$\Bbb E (f(X)) + \Bbb E (f(X_{-})) \ge 2\Bbb E (f(Y)),$$
and the result follows from the property of the function $f$
that $f(X)=f(X_{-}).$ \qed \enddemo

\bigbreak

Let us show that the norm of every subspace of the spaces
$L_r,\ 0<r\le 2$ has all the properties of 
the function $f$ in Theorem 1.

\proclaim{Proposition 2} Let $(\Bbb R^n,\|\cdot\|)$ be a subspace
of $L_r([0,1])$ with $0<r\le 2.$ Then, for any $p\in (-n,0)$
the function $\|x\|^p$ is a positive definite distribution 
on $\Bbb R^n,$ and, therefore, $\Bbb E (\|X\|^p) \ge \Bbb E (\|Y\|^p).$
\endproclaim

\demo{Proof} By the result of Bretagnolle, Dacunha-Castelle 
and Krivine \cite{1}, the
function $\exp(-\|x\|^r)$ is a positive definite function in 
$\Bbb R^n.$ It is easy to see that
$$\|x\|^p = {r\over{\Gamma(-p/r)}}
\int_0^{\infty} |t|^{-1-p} \exp(-|t|^r\|x\|^r)\ dt.$$
The integral in the right-hand side converges because $p<0.$ 
Also that integral represents a positive definite function of $x,$
since $\exp(-|t|^r\|x\|^r)$ is a positive definite function of $x$ 
for every $t>0.$ The inequality for the expectations follows
from Theorem 1. \qed \enddemo 

\bigbreak

We are going to show now that the number of spaces for which
$\|x\|^p$ is a positive definite distribution in $\Bbb R^n$
becomes quite large when $p\to -\infty.$ For example, the norm 
of every $n$-dimensional Banach space has this property if
$-n<p<-n+1.$ 

\proclaim{Proposition 3} Let $-n<p<-n+1.$  Then every even, continuous,
positive, homogeneous of the order $p$ function $f$ on 
$\Bbb R^n\setminus \{0\}$
is a positive definite distribution in $\Bbb R^n.$ 
\endproclaim 

\demo{Proof} Let $\phi$ be any non-negative function from
$\Cal S(\Bbb R^n).$ Writing the integral in the spherical 
coordinates we get
$$(\hat{f},\phi) = \int_{\Bbb R^n} f(x)\hat{\phi}\ dx =
(1/2)\int_S \int_{\Bbb R} f(r\theta) 
|r|^{n-1} \hat{\phi}(r\theta)\ dr\ d\theta =$$
$$(1/2)\int_S f(\theta)\ d\theta 
\Big(\int_{\Bbb R} |r|^{n+p-1}\hat{\phi}(r\theta)\ dr\Big).\tag{13}$$
The integral over $\Bbb R$ in (13) converges because $n+p-1\in (-1,0).$
The Fourier transform of the distribution $|r|^{n+p-1}$ is equal to
$(|r|^{n+p-1})^{\wedge}(t) = c_{n+p-1} |t|^{-n-p},\ t\in \Bbb R$ where 
$c_{n+p-1}=2^{n+p}\sqrt{\pi}\Gamma((n+p)/2)/\Gamma((-n-p+1)/2)$
is a positive constant. 
On the other hand, by the connection between the Fourier transform
and the Radon transform, the function $r\to \hat{\phi}(r\theta)$
is the one-dimensional Fourier transform of the function
$t\to \int_{(x,\theta)=t} \phi(x)\ dx$ ( the latter function is
the Radon transform of $\phi$ in the direction of $\theta.)$
Therefore, for any non-negative function $\phi\in \Cal S(\Bbb R^n),$
switching to the Fourier transforms we get
$$\int_{\Bbb R} |r|^{n+p-1}\hat{\phi}(r\theta)\ dr =
\big(|r|^{n+p-1},\hat{\phi}(r\theta)\big) = $$
$$c_{n+p-1} \big(|t|^{-n-p}, \int_{(x,\theta)=t} \phi(x)\big)=
c_{n+p-1} \int_{\Bbb R} |t|^{-n-p}
\big(\int_{(x,\theta)=t} \phi(x)\ dx\big)\ dt \ge 0,$$
where the last integral converges since $-n-p\in (-1,0).$
We conclude that the integral (13) is non-negative, which means
that $\hat{f}$ is a positive distribution on $\Bbb R^n.$
\qed \enddemo

An immediate consequence of Theorem 1 and Proposition 3 is the
following

\proclaim{Corollary 2} Let $-n<p<-n+1$ and $f$ be any even, continuous,
positive, homogeneous of the order $p$ function in 
$\Bbb R^n\setminus \{0\}$
such that $f(u,v)=f(u,-v)$ 
for every $u\in \Bbb R^k,\ v\in \Bbb R^{n-k}.$ Then
$\Bbb E (f(X)) \ge \Bbb E (f(Y)).$ \endproclaim

Putting $f(x)= \max_{i=1,...,n} |x_i|^p,\ p\in (-n,-n+1)$  in
Corollary 2 we get the inequality (1):

\proclaim{Corollary 3} For any $p\in (-n,-n+1)$ and $q,k, X,Y$ 
as in the Introduction, we have
$$\Bbb E\ (\max_{i=1,...,n} |X_i|)^p \ge 
\Bbb E\ (\max_{i=1,...,n} |Y_i|)^p.$$
\endproclaim

\subheading{Acknowledgements} I would like to thank 
G. Schechtman, T. Schlumprecht and J. Zinn for very useful
discussions.

\enddocument